\documentclass[aos,MSNbibl,seceqn,dvips]{arximspdf}
\usepackage{graphicx}
\doi{10.1214/14-AOS1257} 
\volume{42}
\issue{6}
\pubyear{2014}
\firstpage{2413}
\lastpage{2440}
\docsubty{FLA}

\makeatletter
\newproclaim{remark}{Remark}
\newcommand{\V}[1]{\mathbf{#1}}
\def\cal{\mathcal}
\newcommand{\eqref}[1]{(\ref{#1})}
\newtheorem{thmm}{Theorem}

\newproclaim{defn}{Definition}
\newtheorem{lemma}{Lemma}

\newcommand{\on}{\Delta}
\newcommand{\n}{n} 
\newcommand{\m}{m} 
\newcommand{\rk}{r} 
\newcommand{\Mnn}{M_{\n\times\n}}
\newcommand{\Mmn}{M_{\m\times\n}}
\newcommand{\Mmm}{M_{\m\times\m}}
\newcommand{\MSE}{\mathbf{M}}
\newcommand{\mmx}{\mathcal{M}}
\newcommand{\X}{\mathbf{X}}
\newcommand{\Xmn}{\mathbf{X}_{\m,\n}}
\newcommand{\z}{\zeta}
\newcommand{\cM}{\mmx}
\newcommand{\bR}{\mathbf{R}}
\newcommand{\bX}{\mathbf{X}}
\newcommand{\cA}{{\cal A}}
\newcommand{\goto}{\rightarrow}
\newcommand{\by}{\V{y}}
\newcommand{\bz}{\V{z}}
\newcommand{\eps}{\varepsilon}
\newcommand{\sgn}{\operatorname{sign}}
\newcommand{\cD}{{\cal D}}
\newcommand{\peeOne}{P_1}
\newcommand{\R}{\mathbb{R}}
\newcommand{\E}{\mathbb{E}}
\newcommand{\Nc}{\mathcal{N}}

\makeatother

\begin{document}
\begin{frontmatter}

\title{Minimax risk of matrix denoising by singular value~thresholding}
\runtitle{Minimax risk of matrix denoising}

\begin{aug}
\author{\fnms{David}~\snm{Donoho}\thanksref{NSF}\ead[label=e1]{donoho@stanford.edu}}
\and
\author{\fnms{Matan}~\snm{Gavish}\corref{}\thanksref{NSF,SGF}\ead[label=e2]{gavish@stanford.edu}}
\runauthor{D. Donoho and M. Gavish}
\thankstext{NSF}{Supported in part by NSF Grant DMS-09-06812 (ARRA).}
\thankstext{SGF}{Supported in part by a William R. and Sara Hart Kimball
Stanford Graduate Fellowship and a~Technion EE Sohnis Promising
Scientist Award.}

\affiliation{Stanford University}

\address{Department of Statistics\\
Stanford University\\
Sequoia Hall, 390 Serra Mall\\
Stanford, California 94305-4065\\
USA\\
\printead{e1}\\
\phantom{E-mail:\ }\printead*{e2}}
\end{aug}

\received{\smonth{4} \syear{2013}}
\revised{\smonth{7} \syear{2014}}

%
\begin{abstract}

An unknown $m$ by $n$ matrix $X_0$ is to be estimated from noisy measurements
$Y = X_0 + Z$, where the noise matrix $Z$ has i.i.d. Gaussian entries.
A~popular matrix denoising scheme solves the nuclear norm penalization problem
$\operatorname{min}_X \| Y - X \|_F^2/2 + \lambda\|X\|_* $, where $ \|X\|_*$ denotes
the nuclear norm (sum of singular values). This is the analog, for matrices,
of $\ell_1$ penalization in the vector case. It has been empirically observed
that if $X_0$ has low rank, it may be recovered quite accurately from the
noisy measurement $Y$.

In a proportional growth framework where the rank $r_n$, number of rows $m_n$
and number of columns $n$ all tend to $\infty$ proportionally to each
other ($
r_n/m_n \goto\rho$, $m_n/n \goto\beta$), we evaluate the asymptotic minimax
MSE $\cM( \rho, \beta) = \lim_{m_n,n \goto\infty} \inf_\lambda
\sup_{\operatorname{rank}(X) \leq r_n} \operatorname{MSE}(X_0,\hat{X}_\lambda)$.\break Our formulas involve
incomplete moments of the quarter- and semi-circle laws ($\beta= 1$, square
case) and the Mar\v{c}enko--Pastur law ($\beta< 1$, nonsquare case). For
finite $m$ and $n$, we show that MSE increases as the nonzero singular values
of $X_0$ grow larger. As a result, the finite-$n$ worst-case MSE, a quantity
which can be evaluated numerically, is achieved when the signal $X_0$ is
``infinitely strong.''

The nuclear norm penalization problem is solved by applying soft thresholding
to the singular values of $Y$. We also derive the minimax threshold, namely
the value $\lambda^*(\rho)$, which is the optimal place to threshold the
singular values.

All these results are obtained for general (nonsquare, nonsymmetric) real
matrices. Comparable results are obtained for square symmetric
nonnegative-definite matrices.
\end{abstract}

%
\begin{keyword}[class=AMS]
\kwd[Primary ]{62C20}
\kwd{62H25}
\kwd[; secondary ]{90C25}
\kwd{90C22}
\end{keyword}

\begin{keyword}
\kwd{Matrix denoising}
\kwd{nuclear norm minimization}
\kwd{singular value thresholding}
\kwd{optimal threshold}
\kwd{Stein unbiased risk estimate}
\kwd{monotonicity of power functions of multivariate tests}
\kwd{matrix completion from Gaussian measurements}
\kwd{phase transition}
\end{keyword}
\end{frontmatter}


\section{Introduction} \label{intro:sec}

Suppose we observe a single noisy matrix $Y$, generated by adding noise
$Z$ to
an unknown matrix $X_0$, so that $Y = X_0 + Z$, where $Z$ is a noise matrix.
We wish to recover the matrix $X_0$ with some bound on the mean squared error
(MSE). This is hopeless when $X_0$ is a completely general matrix,
and the
noise $Z$ is arbitrary; but when $X_0$ happens to be of relatively
low rank,
and the noise matrix is i.i.d. standard Gaussian, one can indeed
guarantee quantitatively accurate recovery. This paper provides explicit
formulas for the best possible guarantees obtainable by a popular,
computationally practical procedure.

Specifically, let $Y$, $X_0$ and $Z$ be $\m$-by-$\n$ real matrices (a
set we denote by
$\Mmn$), and suppose that
$Z$ has i.i.d. entries, $Z_{i,j}\sim\mathcal{N}(0,1)$. Consider the following
nuclear-norm penalization (NNP) problem:
%
\begin{equation}
\label{svt:eq} (\mathit{NNP}) \qquad\hat{X}_\lambda=
\mathop{\operatorname{argmin}}_{X\in\Mmn}
\frac{1}{2}\Vert Y-X \Vert_F^2 + \lambda \Vert X
\Vert_*,
\end{equation}
where $\Vert X \Vert_*$ denotes the sum of singular values of $X\in
\Mmn$,
also known
as the nuclear norm, $\Vert\cdot \Vert_F$ denotes square root of the sum
of squared
matrix entries, also known as the Frobenius norm and $\lambda> 0$ is a penalty
factor. A solution to (NNP) is efficiently computable by modern convex
optimization software \cite{Grant2010}; it shrinks away from $Y$ in the
direction of smaller nuclear norm.

Measure performance (risk) by mean-squared error (MSE). When the
unknown $X_0$
is of known rank $r$ and belongs to a matrix class $\Xmn\subset\Mmn
$, the
minimax MSE of NNP is
%
\begin{equation}
\label{svt-mmx:eq} \mmx_{m,\n}(\rk|\X) = \inf_{\lambda}
\mathop{\sup_{X_0\in\Xmn}}_{\operatorname{rank}(X_0)\leq\rk} \frac{1}{\m\n}
\E_{X_0}\bigl\Vert\hat{X}_\lambda (X_0+Z
)-X_0 \bigr\Vert_F^2,
\end{equation}
namely the worst-case risk of $\hat{X}_{\lambda_*}$, where $\lambda
_*$ is the
threshold for which this worst-case risk is the smallest possible. Here,
$\E_{X_0}$ denotes expectation with respect to the random noise matrix $Z$,
conditional on a given value of the signal matrix $X_0$, and
$\hat{X}_\lambda(X_0+Z)$ denotes the denoiser $\hat{X}_\lambda$
acting on the matrix $X_0+Z$. Note that the symbol $\mathbf{X}$
denotes a matrix class, not a particular matrix.
For square
matrices, $\m=\n$, we write $\mmx_\n(\rk|\X)$ instead of $ \mmx
_{\n,\n}(\rk|\X)
$. In a very clear sense $ \mmx_{m,\n}(\rk|\X) $ gives the best possible
guarantee for the MSE of NNP, based solely on the rank and problem
size, and not
on other properties of the matrix $X_0$.

\subsection{Minimax MSE evaluation}
In this paper, we calculate the minimax MSE $ \mmx_{m,\n}(\rk|\X)$
for two matrix classes $\X$:
\begin{longlist}[(1)]
\item[(1)] \textit{General matrices}: $\X=\mathrm{Mat}_{\m,\n}$:
The signal $X_0$ is a real matrix $X_0\in\Mmn$ ($\m\leq\n$).

\item[(2)] \textit{Symmetric matrices}: $\X=\mathrm{Sym}_{\n}$: The signal $X_0$
is a real, symmetric
positive semidefinite matrix, a set we denote by $S^\n_+ \subset\Mnn$.
\end{longlist}

In both cases, the asymptotic MSE (AMSE) in the ``large $\n$'' asymptotic
setting admits considerably simpler and more accessible formulas than the
minimax MSE for finite $\n$. So in addition to the finite-$\n$
minimax MSE, we
study the asymptotic setting where a sequence of problem size triplets
$(r_n,m_n,n)$ is indexed by $n \goto\infty$, and where, along this
sequence $m/n
\goto\beta\in(0,1)$ and $r/m \goto\rho\in(0,1)$.\vadjust{\goodbreak} We think of
$\beta$ as
the matrix shape parameter; $\beta=1$ corresponds to a square matrix, and
$\beta< 1$ to a matrix wider than it is tall. We think of $\rho$ as the
fractional rank parameter, with $\rho\approx0$ implying low rank
relative to
matrix size. Using these notions we can define the asymptotic minimax
MSE (AMSE)
\[
\cM(\rho,\beta|\X) = \lim_{n\to\infty} \mmx_{m_n , \n} (
r_n |\X).
\]

We obtain explicit formulas for the asymptotic minimax MSE in terms of
incomplete moments of classical probability distributions: the quarter-circle
and semi-circle laws (square case $\beta=1$) and the Mar\v{c}enko--Pastur
distribution (nonsquare case $\beta< 1$). Figures~\ref{amse_mat:fig} and \ref{amse_mat_sym:fig} show how
the AMSE
depends on the matrix class $\bX$, the rank fraction $\rho$ and the
shape factor
$\beta$. We also give explicit formulas for the optimal regularization
parameter $\lambda_*$, also as a function of $\rho$; see Figures~\ref{mmx_lambda_mat:fig}
and \ref{mmx_lambda_mat_sym:fig}.

These minimax MSE results constitute best possible guarantees, in the
sense that for
the procedure in question, the MSE is actually attained at some rank $r$
matrix, so that no better guarantee is possible for the given tuning parameter
$\lambda_*$; but also, no other tuning parameter offers a better such
guarantee.

\subsection{Motivations}

We see four reasons to develop these bounds.

\subsubsection{Applications} Several important problems in modern
signal and
image processing, in network data analysis and in computational biology
can be cast as recovery of low-rank matrices from noisy data, and
nuclear norm minimization has become a
popular strategy in many cases; see, for example, \cite{Shabalin2010,Cand2012} and
references therein. Our results provide sharp limits on what such
procedures can
hope to achieve, and validate rigorously the idea that \textit{low rank
alone} is
enough to provide some level of performance guarantee; in fact, they precisely
quantify the best possible guarantee.

\subsubsection{Limits on possible improvements} One might wonder
whether some other procedure offers even better
guarantees than NNP. Consider then
the minimax risk \textit{over all procedures}, defined by
%
\begin{equation}
\label{mmx:eq} \mmx_{m,\n}^*(\rk|\X) = \inf_{\hat{X}}
\mathop{\sup_{X_0\in\Xmn}}_{\operatorname{rank}(X_0)\leq\rk} \frac{1}{\m\n}
\E_{X_0}\bigl\Vert\hat{X} (X_0+Z )-X_0 \bigr\Vert
_F^2 ,
\end{equation}
where $\hat{X} = \hat{X}(Y)$ is some measurable function of
the observations, and its corresponding minimax AMSE
\[
\cM^*(\rho,\beta|\X) = \lim_{n\to\infty} \mmx^*_{m_n , \n} (
r_n |\X) ,
\]
where the sequences $m_n$ and $r_n$ are as above.
Here one wants to find the best possible
procedure, without regard to efficient computation.
We also prove a lower bound on the minimax MSE
over all procedures, and provide an asymptotic evaluation
\[
\cM^*(\rho,\beta| \bX) \geq\cM^-(\rho,\beta) \equiv\rho+ \beta\rho- \beta
\rho^2.
\]
In the square case ($\beta=1$), this simplifies to
$\cM^*(\rho| \bX) \geq\cM^-(\rho) \equiv\rho(2-\rho)$.
The NNP-minimax MSE
is by definition larger than the minimax MSE, $\cM(\rho,\beta|\X)
\geq
\cM^*(\rho,\beta|\X)$.
While there may be
procedures outperforming NNP, the performance improvement
turns out to be limited. Indeed, our formulas show that
\[
\frac{ \cM(\rho,\beta|\X)}{ \cM^-(\rho,\beta)} \leq 2
 \biggl(1 + \frac{\sqrt{\beta}}{1+\beta} \biggr) ,
\]
while
%
\begin{equation}
\label{eq:bndCompare} \lim_{\rho\goto0} \frac{ \cM(\rho,\beta|\X)}{ \cM^-(\rho
,\beta)}= 2 \biggl(1 +
\frac{\sqrt{\beta}}{1+\beta} \biggr).
\end{equation}
For square matrices ($\beta=1$), this simplifies to
%
\begin{equation}
\frac{ \cM(\rho|\X)}{ \cM^-(\rho)} \leq3,\qquad \lim_{\rho\goto0} \frac{ \cM(\rho|\X)}{ \cM^-(\rho)}= 3.
\end{equation}
In words, the potential improvement in minimax AMSE of \textit{any} other
matrix denoising procedure
over NNP is at most a factor of $3$; and if any such improvement were
available, it would only be available in extreme low-rank situations.
Actually obtaining such an improvement in performance guarantees
is an interesting research challenge.

\subsubsection{Parallels in minimax decision theory} The low-rank matrix
denoising problem stands in a line of now-classical problems in minimax
decision theory. Consider the sparse vector denoising problem, where an
unknown vector $x$ of interest yields noisy observations $\by= \mathbf{x}+
\bz$ with
noise $\bz\sim_{\mathrm{i.i.d}.} N(0,1)$; the vector $\mathbf{x}$ is sparsely
nonzero---$
\#\{ i \dvtx x(i) \neq0\} \leq\eps\cdot n$---with $\bz$ and $\mathbf{x}$
independent.
In words, a vector with a fraction $\leq\eps$ of nonzeros is observed with
noise. In this setting, consider the following $\ell_1$-norm penalization
problem:
%
\begin{equation}
\label{soft:eq} (\peeOne)\qquad \hat{\mathbf{x}}_\lambda=
 \mathop{\operatorname{argmin}}_{x\in\bR^n}
\frac{1}{2}\Vert\by-\mathbf{x} \Vert _2^2 + \lambda
\Vert\mathbf{x} \Vert_1.
\end{equation}

The sparse vector denoising problem
exhibits several striking structural resemblances to low-rank matrix denoising:
\begin{itemize}
\item\textit{Thresholding representation}.
For a scalar $y$, define the soft thresholding nonlinearity by
\[
\eta_\lambda(y) = \sgn(y) \cdot\bigl(|y| - \lambda\bigr)_+.
\]
In words, values larger than $\lambda$ are shifted toward zero by
$\lambda$,
while those smaller than $\lambda$ are set to zero. The solution
vector $\hat{x}_\lambda$ of ($\peeOne$)
obeys $(\hat{\mathbf{x}}_\lambda)_i = \eta_\lambda(y_i)$; namely, it
applies $\eta_\lambda$
coordinate wise.
Similarly, the solution of (NNP) applies $\eta_\lambda$
coordinate wise to the singular values of the noisy matrix $Y$.

\begin{remark*}
By this observation, ($\peeOne$) can also be called
``soft thresholding''
or ``soft threshold denoising,'' and in fact, these other terms are the
labels in
common use. Similarly, NNP amounts to ``soft thresholding of singular
values.'' This paper will henceforth use the term \textit{singular value
soft thresholding}
(SVST).
\end{remark*}

\item\textit{Sparsity/low rank analogy.} The objects to be recovered in the
sparse vector denoising problem have sparse entries; those to be
recovered in
the low-rank matrix denoising problem have sparse singular values. Thus the
fractional sparsity parameter $\eps$ is analogous to the fractional rank
parameter $\rho$. It is natural to ask the same questions about
behavior of minimax MSE in
one setting (say, asymptotics as $\rho\goto0$) as in the other
setting ($\eps\goto0$).
In fact, such comparisons turn out to be illuminating.

\item\textit{Structure of the least-favorable estimand.}
Among sparse vectors $x$ of a given fixed sparsity fraction $\eps$,
which of these is the hardest to estimate? This should maximize the
mean-squared error
of soft thresholding, even under the most clever choice of $\lambda$.
This least-favorable configuration is singled out in the minimax AMSE
%
\begin{equation}
\label{eq:sparseMinMax} M_n(\eps) = \inf_{\lambda} \sup
_{\#\{i: x(i) \neq0
\}\leq\eps\cdot n} \frac{1}{n} \E\| \hat{x}_\lambda- x
\|_2^2 .
\end{equation}
In this min/max, the least favorable situation has all its nonzeros, in some
sense, ``at infinity''; that is, all sparse vectors which place large
enough values
on the nonzeros are nearly least favorable, that is, essentially make
the problem
maximally difficult for the estimator, even when it is optimally tuned. In
complete analogy, in low-rank matrix denoising we will see that all low-rank
matrices, which are in an appropriate sense ``sufficiently large,'' are thereby
almost least favorable.

\item\textit{Structure of the minimax smoothing parameter.}
In the sparse vector denoising AMSE (\ref{eq:sparseMinMax}) the
$\lambda= \lambda(\eps)$ achieving the infimum is a type of optimal
regularization parameter, or optimal threshold. It decreases as $\eps$
increases, with $\lambda(\eps) \goto0$ as $\eps\goto1$.
Paralleling this, we show that the low-rank matrix denoising AMSE~(\ref{svt-mmx:eq}) has minimax singular value soft threshold $ \lambda
^*(\rho)$ decreasing as $\rho$ increases,
and $\lambda^*(\rho) \goto0$ as $\rho\goto1$.
\end{itemize}

Despite these similarities, there is one major difference between
sparse vector denoising and low-rank matrix denoising.
In the sparse vector denoising problem, the soft-thresholding minimax MSE
was compared to the minimax MSE over all procedures by Donoho and
Johnstone \cite{Donoho1994}.
Let $M(\eps) = \lim_{n \goto\infty} M_n(\eps)$
denote the soft thresholding AMSE and define the minimax AMSE over all
procedures via
\[
M^*(\eps) = \lim_{n \goto\infty} \inf_{\hat{x}} \sup
_{\#\{i: x(i) \neq0 \}\leq\eps\cdot n} \frac{1}{n}\E\| \hat{x} - x \|_2^2,
\]
where here $\hat{x} = \hat{x}(y)$ denotes \textit{any} procedure which
is measurable in the observations.
In the limit of extreme sparsity, soft thresholding is
\textit{asymptotically minimax}~\cite{Donoho1994},
\[
\frac{M(\eps)}{M^*(\eps)} \goto1 \qquad\mbox{as } \eps\goto0.
\]
Breaking the chain of similarities, we are not able to
show a similar asymptotic minimaxity for
SVST in the low rank matrix denoising problem.
Although equation (\ref{eq:bndCompare}) says that
soft thresholding of singular values is asymptotically
not more than a factor of 3 suboptimal, we doubt that
anything better than a factor of $3$ can be true; specifically,
we conjecture that SVST suffers a \textit{minimaxity gap}. For example, for
$\beta=1$, we conjecture that
\[
\frac{\mmx(\rho|\X)}{\mmx^*(\rho|\X)} \goto3\qquad \mbox{as } \rho \goto0.
\]
We believe that interesting new estimators will be found improving upon
singular value soft thresholding by essentially this factor of $3$.
Namely, there may
be substantially better guarantees to be had under extreme sparsity,
than those which can be offered
by SVST.
Settling the minimaxity gap for SVST seems a challenging new research question.

\subsubsection{Indirect observations}

Evaluating the Minimax MSE of SVST has
an intriguing new motivation \cite{Donoho2011,Donoho2013,OymakHassibi12},
arising from the newly evolving fields of compressed sensing and matrix
completion.

Consider the problem of recovering an unknown matrix $X_0$ from \textit{noiseless}, \textit{indirect} measurements.
Let $\cA\dvtx \R^{m\times n}\to\R^p$ be a linear operator, and consider
observations
\[
y = \cA(X_0).
\]
In words, $y \in\R^p$ contains $p$ linear measurements
of the matrix object $X_0$. See the closely related {\em
trace regression} model \cite{Rohde2011} which also includes
measurement noise. Can we recover $X_0$?
It may seem that $p \geq mn$ measurements are required, and in general
this would be true;
but if $X_0$ happens to be of low rank, and $\cA$ has suitable properties,
we may need substantially fewer measurements.

Consider reconstruction by \textit{nuclear norm minimization},
%
\begin{equation}
\label{NNMC} (P_{\mathrm{nuc}}) \qquad \operatorname{min} \|X\|_*
\qquad\mbox{subject to }
y=\cA(X).
\end{equation}

Recht and co-authors found that when the matrix representing the
operator $\cA$
has i.i.d. $\Nc(0,1)$ entries, and the matrix is of rank $r$,
the matrix $X_0$ is recoverable from $p < nm$ measurements
for certain combinations of $p$ and $r$ \cite{Recht2010}.
The operator $\cA$ offers so-called \textit{Gaussian measurements} when the
representation of the operator as a matrix has i.i.d. Gaussian entries.
Empirical work by Recht, Xu and Hassibi~\cite{Recht2010a,Recht2008}, Fazel,
Parillo and Recht \cite{Recht2010},
Tanner and Wei \cite{Tanner2012} and
Oymak and Hassibi \cite{OymakHassibi11}
documented for Gaussian measurements a \textit{phase transition}
phenomenon, that is,
a fairly sharp transition from success to failure as $r$ increases, for
a given $p$.
Putting $\rho= r/m$ and $\delta= p/(mn)$ it appears that there is
a critical sampling rate $\delta^*(\rho)=\delta^*(\rho;\beta)$,
such that, for $\delta> \delta^*(\rho)$,
NNM is successful for large~$m,n$, while for $\delta< \delta^*(\rho)$,
NNM fails. $\delta^*(\rho)$
provides a sharp ``sampling limit'' for low rank matrices, that is,
a clear statement of how many measurements are needed to recover a low
rank matrix,
by a popular and computationally tractable algorithm.

In very recent work, \cite{Donoho2011,Donoho2013,OymakHassibi12}, it
has been shown empirically
that the precise location of the phase transition \textit{coincides with
the minimax MSE}
%
\begin{equation}
\label{eq:PTFormula} \delta^*(\rho; \beta) = \mmx(\rho,\beta|\X),\qquad \rho\in(0,1) , \beta
\in(0,1);
\end{equation}
a key requirement for discovering and verifying (\ref{eq:PTFormula})
empirically was to obtain an explicit formula for
the right-hand side; that explicit formula is derived and proven in
this paper. Relationship (\ref{eq:PTFormula})
connects two seemingly unrelated problems: matrix denoising from direct
observations and matrix recovery from
incomplete measurements. Both problems are attracting a large and
growing research literature.
Equation (\ref{eq:PTFormula}) demonstrates the importance of minimax
MSE calculations
even in a seemingly unrelated setting where there is no noise and no
statistical decision to be made!

\section{Results} \label{results:sec}

\subsection{Least-favorable matrix}

We start by identifying the least-favorable situation for matrix
denoising by
SVST.

\begin{thmm}[(The worst-case matrix for SVST has its principal subspace ``at
$\infty$'')] \label{lf:thmm}
Define the risk function of a denoiser
$\hat{X}\dvtx \Mmn\to\Mmn$ at $X_0\in\Mmn$ by
%
\begin{equation}
R(\hat{X},X_0):= \frac{1}{\m}\E\biggl\Vert\hat{X}
\biggl(X_0 + \frac
{1}{\sqrt{\n}}Z \biggr) - X_0
\biggr\Vert_F^2. \label{risk:eqn}
\end{equation}
Let $\lambda>0$, $\m\leq\n\in\mathbb{N}$ and $1\leq\rk\leq\m$.
For the worst-case risk of $\hat{X}_\lambda$ on $\m\times\n$
matrices of rank
at most $\rk$,
we have
%
\begin{equation}
\mathop{\sup_{X_0\in\Mmn}}_{\operatorname{rank}(X_0)\leq\rk} R(\hat{X}_\lambda,X_0)
= \lim_{\mu\to\infty} R(\hat{X}_\lambda, \mu C) ,
\end{equation}
where $C\in\Mmn$ is any fixed matrix of rank exactly $\rk$.
\end{thmm}

\subsection{Minimax MSE}

Let $W_i(\m,\n)$ denote the marginal distribution of the $i$th
largest eigenvalue of a standard central Wishart matrix $W_\m(I,\n)$,
namely, the $i$th
largest eigenvalue of the random matrix $\frac{1}{\n}ZZ'$ where $Z\in
\Mmn$ has
i.i.d. $\mathcal{N}(0,1)$ entries.
Define for $\Lambda>0$ and $\alpha\in \{ 1/2,1  \}$
%
\begin{eqnarray}
\label{finite-n-proxy:eq} \MSE_{\n}(\Lambda;\rk,\m,\alpha) &=&
\frac{\rk}{\m}+\frac{\rk}{\n}-\frac{r^2}{\m\n} + \frac{\rk(\n-\rk)}{\m\n}
\Lambda^2
\nonumber
\\[-8pt]
\\[-8pt]
\nonumber
&&{}+ \alpha\frac{(\n-\rk)}{\m\n} \sum_{i=1}^{\m-\rk}w_i
(\Lambda;\m-\rk;\n-\rk ) ,
\end{eqnarray}
where
%
\begin{equation}
\label{wishart-eig-moment:eq} w_i ( \Lambda;\m,\n ) = \int_{\Lambda^2}^\infty
(\sqrt{t}-\Lambda)^2 \,dW_i(\m,\n) (t)
\end{equation}
is a combination of the complementary incomplete moments of standard
central Wishart
eigenvalues
\[
\int_{\Lambda^2}^\infty t^{k/2} \,dW_i(\m,\n)
(t)
\]
for $k=0,1,2$.

\begin{thmm}[(An implicit formula for the finite-$\n$ minimax MSE)] \label{finite-n-mmx:thmm}
The minimax MSE of SVST over $\m$-by-$\n$ matrices of rank at most
$\rk$ is given by
\begin{eqnarray*}
\mmx_\n(\rk,\m|\mathrm{Mat}) &=& \min_{\Lambda\geq0}
\MSE_\n(\Lambda;\rk,\m,1)\qquad \mbox{and}
\\
\mmx_\n(\rk|\mathrm{Sym}) &=& \min_{\Lambda\geq0}
\MSE_\n(\Lambda;\rk ,\n,1/2) ,
\end{eqnarray*}
where the minimum on the right-hand sides is unique.
\end{thmm}

In fact, we will see that $ \MSE_\n(\Lambda;\rk,\m,\alpha)$
is convex in $\Lambda$.
As the densities of the standard central Wishart eigenvalues $W_i(\m
,\n)$ are known
\cite{Zanella2009a}, this makes it possible, in principle, to tabulate
the finite-$n$ minimax
risk.

\subsection{Asymptotic minimax MSE}

A more accessible formula is obtained by calculating the large-$n$ asymptotic
minimax MSE, where $\rk=\rk(\n)$ and $\m=\m(\n)$ both grow
proportionally to
$\n$. Let us write \textit{minimax AMSE} for asymptotic minimax MSE.
For the case $\Xmn=\mathrm{Mat}_{\m,\n}$ we assume a limiting rank fraction
$\rho=
\lim_{\n\to\infty} \rk/\m$ and limiting aspect ratio
$\beta=\lim_{\n\to\infty}\m/\n$ and consider
%
\begin{eqnarray}
\label{asymp-mmx-svt:eq} \mmx(\rho,\beta| \mathrm{Mat}) &=& \lim_{\n\to\infty}
\mmx_\n(\rk,\m|\mathrm{Mat})
\nonumber
\\[-8pt]
\\[-8pt]
\nonumber
&=& \lim_{\n\to\infty}\inf_{\lambda} \mathop{\sup
_{X_0\in M_{\lceil\beta\n\rceil\times n}}}_{\operatorname{rank}(X_0)\leq
\rho\beta\n} \frac{1}{\m\n}\E\Vert
\hat{X}_\lambda-X_0 \Vert_F^2.
\end{eqnarray}
Similarly, for the case $\Xmn=\mathrm{Sym}_{\n}$, we assume a limiting rank fraction
$\rho= \lim_{\n\to\infty} \rk/\n$ and consider
%
\begin{eqnarray}
\mmx(\rho| \mathrm{Sym}) &=& \lim_{n\to\infty} \mmx_\n(\rk|\mathrm{Sym})
\nonumber
\\[-8pt]
\\[-8pt]
\nonumber
&=& \lim_{n\to\infty} \inf_{\lambda} \mathop{\sup
_{
{X_0\in S^\n_+}}}_{\operatorname{rank}(X_0)\leq\rho\n} \frac{1}{\n^2}\E\Vert
\hat{X}_\lambda-X_0 \Vert_F^2.
\end{eqnarray}

The Mar\u{c}enko--Pastur distribution \cite{Marcenko1967}
gives the asymptotic empirical distribution of Wishart eigenvalues. It has
density
%
\begin{equation}
\label{mp:eq} p_\gamma(t) = \frac{1}{2\pi\gamma t} \sqrt{(
\gamma_{+} - t) (t-\gamma_{-}}) \cdot\mathbf{1}_{[\gamma_-,\gamma_+]}(t),
\end{equation}
where $\gamma_{\pm}= ( 1\pm\sqrt{\gamma}  )^2$.
Define the
complementary incomplete moments of the Mar\u{c}enko--Pastur distribution
%
\begin{equation}
\label{MP:eq} P_{\gamma}(x;k)= \int_x^{\gamma_+}
t^{k} p_\gamma(t) \,dt.
\end{equation}
Finally, let
%
\begin{eqnarray}
\label{proxy:eq} &&\mathbf{M}(\Lambda;\rho,\tilde{\rho},\alpha)\nonumber\\
&&\qquad= \rho+ \tilde{
\rho} - \rho\tilde{\rho} + (1-\tilde{\rho})
\nonumber\qquad
\\[-8pt]
\\[-8pt]
\nonumber
&&\qquad\quad{}\times \biggl[ \rho\Lambda^2 \\
&&\hspace*{16pt}\qquad\quad{}+\alpha(1-\rho) \biggl( P_\gamma\bigl( \Lambda^2 ; 1
\bigr) - 2\Lambda P_\gamma\biggl( \Lambda^2 ;
\frac{1}{2}\biggr) + \Lambda^2 P_\gamma\bigl(
\Lambda^2 ; 0\bigr) \biggr) \biggr] ,\nonumber
\end{eqnarray}
with
$\gamma=\gamma(\rho,\tilde{\rho}) = (\tilde{\rho} - \rho\tilde
{\rho}) /
(\rho-\rho\tilde{\rho})$.

\begin{thmm}[(An explicit formula for the minimax AMSE)] \label{asymp-mmx:thmm}
For the minimax AMSE of SVST we have
%
\begin{eqnarray}
\mmx(\rho,\beta| \mathrm{Mat}) &=& \min_{0\leq\Lambda\leq\gamma_+} \mathbf{M}(\Lambda; \rho,
\beta\rho,1), \label{asymp-mmx-mat:eq}
\\
\mmx(\rho| \mathrm{Sym}) &=& \min_{0\leq\Lambda\leq\gamma_+} \mathbf {M}(\Lambda; \rho,
\rho,1/2) , \label{asymp-mmx-sym:eq}
\end{eqnarray}
with $\gamma_+ =  (1+\sqrt{(\beta-\beta\rho)/(1-\beta\rho
)} )^2$,
where the minimum on the right-hand sides is unique.
Moreover, for any $0<\beta\leq1$, the function
$\rho\mapsto\mmx(\rho,\beta|\mathrm{Mat})$
is continuous and increasing on $\rho\in[0,1]$, with
$\mmx(0,\beta|\mathrm{Mat})=0$ and
$\mmx(1,\break \beta|\mathrm{Mat})=1$. The same is true for $\mmx(\rho|\mathrm{Sym})$.
\end{thmm}

\begin{figure}

\includegraphics{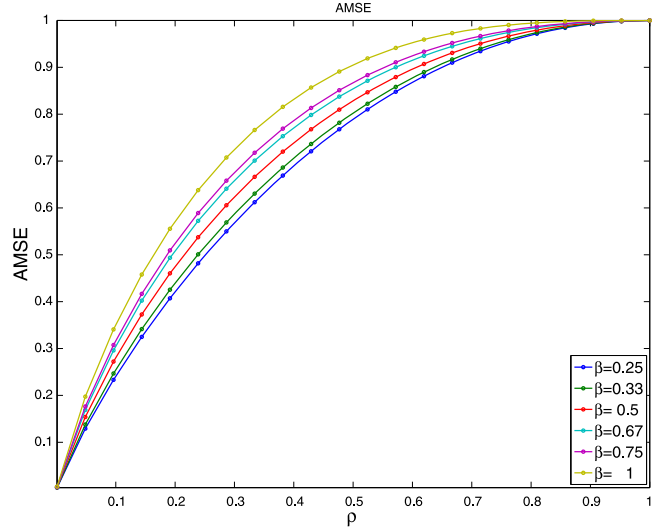}

\caption{The minimax AMSE curves for case $\mathrm{Mat}$,
defined in \protect\eqref{asymp-mmx-mat:eq},
for a few values of $\beta$.}
\label{amse_mat:fig}
\end{figure}

\begin{figure}

\includegraphics{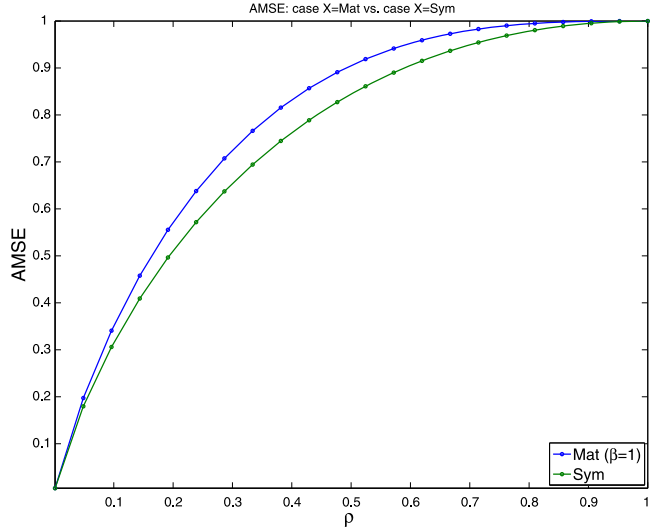}

\caption{The minimax AMSE curves for case $\mathrm{Mat}$
with $\beta=1$ and case $\mathrm{Sym}$.}
\label{amse_mat_sym:fig}
\end{figure}

The curves $\rho\mapsto\mmx(\rho,\beta| \mathrm{Mat})$, for different
values of
$\beta$, are shown in Figure~\ref{amse_mat:fig}. The curves
$\rho\mapsto\mmx(\rho,\beta| \mathrm{Mat})$ and $\rho\mapsto\mmx(\rho
,\beta| \mathrm{Mat})$
are shown in Figure~\ref{amse_mat_sym:fig}.

\subsection{Computing the minimax AMSE}

To compute $\mmx(\rho,\beta|\mathrm{Mat})$ and\break $\mmx(\rho|\mathrm{Sym})$ we need to
minimize \eqref{proxy:eq}.
Define
%
\begin{equation}
\label{proxy-minimizer:eq} \Lambda_*(\rho,\beta,\alpha) =
\mathop{\operatorname{argmin}}_\Lambda
\mathbf{M}(\Lambda;\rho,\tilde{\rho },\alpha).
\end{equation}

\begin{thmm}[(A characterization of the minimax AMSE for general $\beta$)] \label{proxy-minimizer:thmm}
For any $\alpha\in \{ 1/2,1  \}$ and $\beta\in(0,1]$,
the function
$\rho\mapsto\Lambda_*(\rho,\beta,\alpha)$ is decreasing on
$\rho\in[0,1]$ with\vspace*{2pt}
%
\begin{eqnarray}\label{lim_Lambda_down:eq}
\lim_{\rho\to0}\Lambda_*(\rho,\beta,\alpha) &=& \Lambda_*(0,\beta,
\alpha) = 1 + \sqrt{\beta}  \qquad\mbox{and}
\\[2pt]
\lim_{\rho\to1}\Lambda_*(\rho,\beta,\alpha) &=& \Lambda_*(1,\beta,
\alpha) = 0. \label{lim_Lambda_up:eq}\vspace*{2pt}
\end{eqnarray}
For $\rho\in(0,1)$,
the minimizer $\Lambda_*(\rho,\beta,\alpha)$ is the unique root of the
equation in $\Lambda$\vspace*{2pt}
%
\begin{equation}
\label{argmin-Lambda:eq} P_\gamma\biggl(\Lambda^2;\frac{1}{2}
\biggr) - \Lambda\cdot P_\gamma\bigl(\Lambda^2;0\bigr) =
\frac{\Lambda\rho}{\alpha
(1-\rho)} ,\vspace*{2pt}
\end{equation}
where the left-hand side of \eqref{argmin-Lambda:eq} is a
decreasing function of $\Lambda$.\vspace*{3pt}
\end{thmm}

The minimizer $\Lambda_*(\rho,\beta,\alpha)$ can therefore be determined
numerically by binary search.
[In fact, we will see that $\Lambda_*$ is the unique minimizer of the convex
function $\Lambda\mapsto\MSE(\Lambda; \rho,\tilde{\rho},\alpha)$.]
Evaluating $\mmx(\rho,\beta|\mathrm{Mat})$ and
$\mmx(\rho|\mathrm{Sym})$ to precision $\epsilon$
thus requires $O(\log(1/\epsilon))$ evaluations of the complementary
incomplete Mar\u{c}enko--Pastur moments \eqref{MP:eq}.\vspace*{1pt}

For square matrices ($\beta=1$), this computation turns out to be even simpler,
and only requires
evaluation of elementary trigonometric functions.\vspace*{3pt}
%
\begin{thmm}[(A characterization of the minimax AMSE for $\beta=1$)] \label{proxy-sq:thmm}
We have\vspace*{1pt}
%
\begin{eqnarray}
\label{proxy-sq:eq} \mathbf{M}(\Lambda;\rho,\rho,\alpha)&=&\rho(2-\rho)
\nonumber
\\[-4pt]
\\[-12pt]
\nonumber
&&{} + (1-\rho)
\bigl[ \rho\Lambda^2
+ \alpha(1-\rho) \bigl( Q_2 ( \Lambda ) - 2\lambda Q_1
( \Lambda ) + \Lambda^2 Q_0 ( \Lambda ) \bigr) \bigr] ,\nonumber\vspace*{1pt}
\end{eqnarray}
where\vspace*{1pt}
%
%
\begin{eqnarray}\label{quarter-circle-moments0:eq}
\label{Q0:eq} Q_{0}(x) &=& \frac{1}{\pi} \int_x^2
\sqrt{4-t^2} \,dt
\nonumber
\\[-7pt]
\\[-7pt]
\nonumber
& =& 1 - \frac{x}{2\pi}\sqrt{4-x^2} -
\frac{2}{\pi} a\operatorname{tan}\biggl(\frac{x}{\sqrt
{4-x^2}} \biggr),
\\[1pt]
Q_{1}(x) &=& \frac{1}{\pi}\int_x^2 t
\sqrt{4-t^2} \,dt = \frac{1}{3\pi}\bigl(4 - x^2
\bigr)^{3/2}, \label{quarter-circle-moments1:eq}
\\[1pt]
\label{Q2:eq} Q_{2}(x) &=& \frac{1}{\pi}\int_x^2
t^2 \sqrt{4-t^2} \,dt
\nonumber
\\[-7pt]
\\[-7pt]
\nonumber
& =& 1 - \frac{1}{4\pi}x
\sqrt{4-x^2}\bigl(x^2-2\bigr) - \frac{2}{\pi}a\operatorname{sin}
\biggl(\frac{x}{2}\biggr) \label{quarter-circle-moments2:eq}\vadjust{\goodbreak}
\end{eqnarray}
are the complementary incomplete moments of the quarter circle law.
Moreover, for $\alpha\in \{ 1/2,1  \}$
%
\begin{equation}
\Lambda_*(\rho,\rho,\alpha) = 2 \cdot \sin \bigl( \theta_\alpha(\rho)
\bigr) ,
\end{equation}
where $\theta_\alpha(\rho)\in[0,\pi/2]$ is the unique solution to
the transcendental equation
%
\begin{equation}
\label{argmin-Lambda-sq:eq} \theta+ \operatorname{cot}(\theta)\cdot \biggl( 1-\frac{1}{3}\operatorname{cos}^2(
\theta) \biggr)= \frac{\pi(1+\alpha^{-1}\rho-\rho)}{2(1-\rho)}.
\end{equation}
The left-hand side of \eqref{argmin-Lambda-sq:eq} is
a decreasing function of $\theta$.
\end{thmm}

In \cite{pt-mmx-runmycode} we make available a Matlab script, and a web-based
calculator for evaluating
$\mmx(\rho,\beta|\mathrm{Mat})$ and $\mmx(\rho|\mathrm{Sym})$. The implementation
provided employs
binary search to solve
\eqref{argmin-Lambda:eq} [or \eqref{argmin-Lambda-sq:eq} for $\beta=1$]
and then feeds the minimizer $\Lambda_*$ into~\eqref{proxy:eq} [or
into \eqref{proxy-sq:eq} for
$\beta=1$].

\subsection{Asymptotically optimal tuning for the SVST threshold
\texorpdfstring{$\lambda$}{lambda}}

The crucial functional $\Lambda_*$, defined in \eqref{proxy-minimizer:eq},
can now be explained as the optimal (minimax)
threshold of SVST in a special system of units.
Let $\lambda_*(\m,\n,\rk|\X)$ denote the minimax tuning threshold, namely
\[
\lambda_*(\m,\n,\rk|\X) = \mathop{\operatorname{argmin}}_{\lambda} \mathop{\sup
_{
{X_0\in\Xmn}}}_{\operatorname{rank}(X_0)\leq\rk} \frac{1}{\m\n}\E_{X_0}\bigl\Vert
\hat{X}_\lambda (X_0+Z )-X_0
\bigr\Vert_F^2.
\]
%

\begin{thmm}[(Asymptotic minimax tuning of SVST)] \label{asymp-tuning-mmx:thmm}
Consider again a sequence $\n\mapsto(\m(\n),\rk(\n))$ with a
limiting rank fraction $\rho=
\lim_{\n\to\infty} \rk/\m$ and a limiting aspect ratio
$\beta=\lim_{\n\to\infty}\m/\n$.
For the asymptotic minimax tuning threshold we have
\begin{eqnarray*}
\lim_{\n\to\infty}\frac{1}{\sqrt{\n}} \lambda_*(\m,\n,\rk |\mathrm{Mat}) &=&
\sqrt{(1-\beta\rho)}\cdot \Lambda_*(\rho,\beta,1) \quad\mbox{and}
\\
\lim_{\n\to\infty}\frac{1}{\sqrt{\n}} \lambda_*(\n,\rk|\mathrm{Sym}) &=&
\sqrt{(1-\rho)}\cdot \Lambda_*(\rho,1,1/2).
\end{eqnarray*}
\end{thmm}

\begin{figure}

\includegraphics{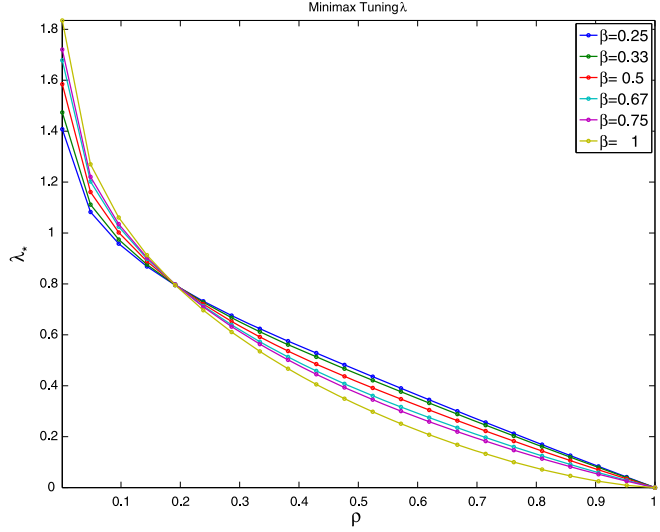}
\caption{(Nonsquare cases.) The scaled asymptotic minimax tuning
threshold for SVST,
$\rho\mapsto\lim_{\n\to\infty}\lambda_*(\m,\n,\rk| \mathrm{Mat})/\sqrt
{\n}$,
when $\m/\n\to\beta$ and $\rk/\m\to\rho$,
for a few values of $\beta$.}
\label{mmx_lambda_mat:fig}
\end{figure}

\begin{figure}
\includegraphics{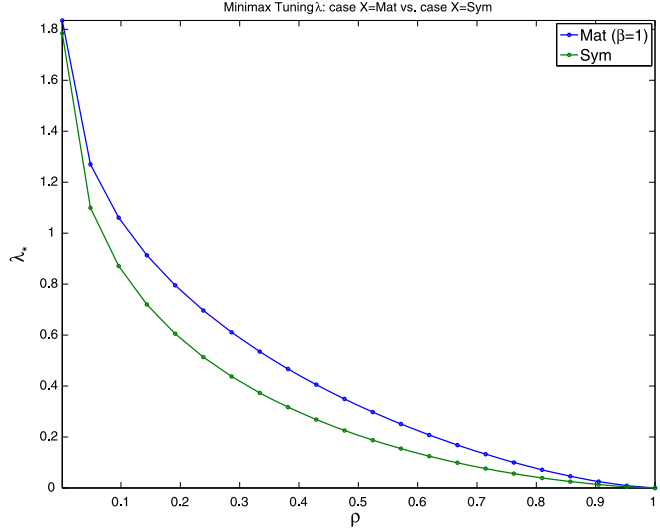}
\caption{(Square case.) The scaled asymptotic minimax tuning
threshold for SVST,
$\rho\mapsto\lim_{\n\to\infty}\lambda_*(\n,\n,\rk| \mathrm{Mat})/\sqrt
{\n}$
and
$\rho\mapsto\lim_{\n\to\infty}\lambda_*(\n,\rk| \mathrm{Sym})/\sqrt{\n}$,
when $\rk/\m\to\rho$.}
\label{mmx_lambda_mat_sym:fig}
\end{figure}

The curves $\rho\mapsto\lim_{\n\to\infty}\lambda_*(\m,\n,\rk|
\mathrm{Mat})/\sqrt{\n}$,
namely the scaled asymptotic minimax tuning threshold for SVST,
are shown in Figure~\ref{mmx_lambda_mat:fig}
for different values of $\beta$.
The curves
$\rho\mapsto\lim_{\n\to\infty}\lambda_*(\n,\n,\rk| \mathrm{Mat})/\sqrt
{\n}$
and
$\rho\mapsto\lim_{\n\to\infty}\lambda_*(\n,\rk| \mathrm{Sym})/\sqrt{\n}$
are shown in Figure~\ref{mmx_lambda_mat_sym:fig}.

\subsection{Parametric representation of the minimax AMSE for square matrices}

For square matrices ($\rho=\tilde{\rho}$, $\beta=1$) the minimax
curves $\mmx(\rho,1|\mathrm{Mat})$
and $\mmx(\rho|\mathrm{Sym})$ admit a parametric representation in the $(\rho
,\mmx)$
plane using elementary trigonometric functions.

\begin{thmm}[(Parametric representation of the minimax AMSE curve for \mbox{$\beta=1$})] \label{parametric:thmm}
As $\theta$ ranges over
$(0,\pi/2)$,
\begin{eqnarray*}
\rho(\theta) &=& 1 - \frac{\pi/2}{
\theta+(\cot(\theta)\cdot(1-({1}/{3})\operatorname{cos}^2(\theta)))},
\\
\mmx(\theta) &=& 2\rho(\theta) - \rho^2(\theta) + 4\rho(\theta)
\bigl(1-\rho(\theta)\bigr)\sin^2(\theta)
\\
&&{}+ \frac{4}{\pi}(1-\rho)^2 \biggl[ (\pi-2\theta) \biggl(
\frac{5}{4} - \operatorname{cos}(\theta)^2\biggr) + \frac{\sin(2\theta
)}{12}
\bigl(\operatorname{cos}(2\theta)-14\bigr) \biggr]
\end{eqnarray*}
is a parametric representation of $\rho\mapsto\mmx(\rho,\rho
|\mathrm{Mat})$, and similarly
\begin{eqnarray*}
\rho(\theta) &=& 1 - \frac{
\theta+(\cot(\theta)\cdot(1-({1}/{3})\operatorname{cos}^2(\theta))) - \pi/2
}{
\theta+(\cot(\theta)\cdot(1-({1}/{3})\operatorname{cos}^2(\theta))) + \pi/2
},
\\
\mmx(\theta) &=& 2\rho(\theta) - \rho^2(\theta) + 4\rho(\theta)
\bigl(1-\rho(\theta)\bigr)\sin^2(\theta)
\\
&&{}+ \frac{2}{\pi}(1-\rho)^2 \biggl[ (\pi-2\theta) \biggl(
\frac{5}{4} - \operatorname{cos}(\theta)^2\biggr) + \frac{\sin(2\theta
)}{12}
\bigl(\operatorname{cos}(2\theta)-14\bigr) \biggr]
\end{eqnarray*}
is a parametric representation of $\rho\mapsto\mmx(\rho|\mathrm{Sym})$.
\end{thmm}

\subsection{Minimax AMSE in the low-rank limit \texorpdfstring{$\rho\approx0$}
{rho approx 0}}

\begin{thmm}[(Minimax AMSE to first order in $\rho$ near $\rho=0$)] \label{small-rho-mmx:thmm}
For the behavior of the minimax curves near $\rho=0$, we have
\[
\mmx(\rho,\beta|\mathrm{Mat}) = 2 ( 1+\sqrt{\beta}+\beta )\cdot\rho+ o(\rho)
\]
and in particular
\[
\mmx(\rho,1|\mathrm{Mat}) = 6\rho+ o(\rho).
\]
Moreover,
\[
\mmx(\rho|\mathrm{Sym}) = 6\rho+ o(\rho).
\]
\end{thmm}

\begin{figure}

\includegraphics{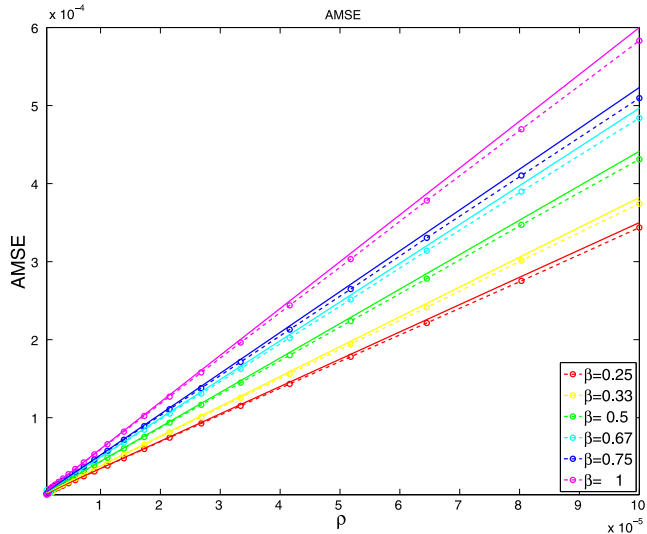}

\caption{The minimax AMSE curves $\rho\mapsto\mmx(\rho,\beta|\mathrm{Mat})$ for
small values of $\rho$
(dashed lines)
and the corresponding approximation slopes $2(1+\sqrt{\beta}+\beta)$ (solid
lines).}
\label{amse_small_rho:fig}
\end{figure}

The minimax AMSE curves $\rho\mapsto\mmx(\rho,\beta|\mathrm{Mat})$ for
small values of
$\rho$, and the corresponding approximation slopes $2(1+\sqrt{\beta
}+\beta)$ are
shown in Figure~\ref{amse_small_rho:fig} for several values of $\beta
$. We find
it surprising that asymptotically, \textit{symmetric positive definite
matrices are
no easier to recover than general square matrices}. This phenomenon is also
seen in the case of sparse vector denoising, where in the limit of extreme
sparsity, the nonnegativity of the nonzeros does not allow one to
reduce the
minimax MSE.\setcounter{footnote}{2}\footnote{Compare results in \cite{Donoho1994} with
\cite{Donoho1992}. To be clear, in both matrix denoising and vector denoising,
there is an MSE advantage for each fixed positive rank fraction/sparsity
fraction. It is just that the benefit goes away as either fraction tends to
$0$.}
We note that this first-order AMSE near $\rho=0$ agrees with a different
asymptotic model for minimax MSE of SVST over large low-rank matrices
\cite{Donoho2013b}. There, the asymptotic prediction for AMSE
near $\rho=0$ is found to be in agreement with the empirical
finite-$n$ MSE.

\subsection{AMSE vs. the asymptotic global minimax MSE}

In \eqref{mmx:eq} we have introduced global minimax MSE
$\mmx^*_{\m,\n}(\rk|\bX)$,
namely the minimax risk over \textit{all} measurable denoisers
$\hat{X}\dvtx \Mmn\to\Mmn$.
To define the large-$n$ asymptotic global minimax MSE analogous to
\eqref{asymp-mmx-svt:eq},
consider sequences where $\rk=\rk(\n)$ and $\m=\m(\n)$
both grow proportionally to
$\n$, such that both limits $\rho=
\lim_{\n\to\infty} \rk/\m$ and
$\beta=\lim_{\n\to\infty}\m/\n$ exist.
Define the asymptotic global minimax MSE
%
\begin{equation}
\mmx^*(\rho,\beta|\bX) = \lim_{\n\to\infty}\mmx^*_{\m,\n
}(\rk|
\bX) .
\end{equation}

\begin{thmm} \label{global-mmx:thmm}
\textup{(1)} For the global minimax MSE we have
%
\begin{equation}
\mmx^*_{\m,\n}(\rk|\bX) \geq \frac{\rk}{\m}+\frac{\rk}{\n}-
\frac{\rk^2+\rk}{\m\n}
\end{equation}
for case $\mathrm{Mat}$, and if $\m=\n$, for case $\mathrm{Sym}$.\vspace*{-6pt}
\begin{longlist}[(3)]
\item[(2)]
For the asymptotic global minimax MSE we have
%
\begin{equation}
\mmx^*(\rho,\beta|\bX) \geq \rho+ \tilde{\rho} - \rho\tilde{\rho}
\end{equation}
for case $\mathrm{Mat}$, and if $\beta=1$, for case $\mathrm{Sym}$. Here
$\tilde{\rho} = \beta\rho$.
\item[(3)]
Let
%
\begin{equation}
\mmx^{-}(\rho,\beta) = \rho+ \tilde{\rho} - \rho\tilde{\rho}
\end{equation}
denote our lower bound on asymptotic global minimax MSE.
Then
%
\begin{equation}
\frac{\mmx(\rho,\beta|\bX)}{\mmx^{-}(\rho,\beta)} \leq 2
\biggl(1 + \frac{\sqrt{\beta}}{1+\beta} \biggr)
\end{equation}
and
%
\begin{equation}
\lim_{\rho\to0} \frac{\mmx(\rho,\beta|\bX)}{\mmx^{-}(\rho,\beta)} = 2 \biggl(1 + \frac{\sqrt{\beta}}{1+\beta}
\biggr).
\end{equation}
\end{longlist}
\end{thmm}

\subsection{Outline of this paper}

The body of the paper proves the above results. Section~\ref{prelim:sec}
introduces notation, and proves auxiliary lemmas. In Section~\ref{lf:sec} we
characterize the worst-case MSE of SVST for matrices of a fixed size (Theorem~\ref{lf:thmm}). In Section~\ref{finite-n-proxy:sec} we derive formula \eqref{finite-n-proxy:eq}
for the
worst-case MSE, and prove Theorem~\ref{finite-n-mmx:thmm}. In Section~\ref{asymp-proxy:sec} we pass to the large-$\n$ limit and derive formula
\eqref{proxy:eq}, which provides the worst-case asymptotic MSE in the
large-$\n$
limit (Theorem~\ref{asymp-mmx:thmm}).
In Section~\ref{asymp-minimizer:sec} we investigate the minimizer of the asymptotic
worst-case MSE function, and its minimum, namely the minimax AMSE, and prove
Theorem~\ref{proxy-minimizer:thmm}.
In Section~\ref{global:sec} we extend our
scope from SVST denoisers to all denoisers, investigate the global
minimax MSE
and prove Theorem~\ref{global-mmx:thmm}.
In the interest of space, Theorems \ref{proxy-sq:thmm},
\ref{asymp-tuning-mmx:thmm}
\ref{parametric:thmm} and
\ref{small-rho-mmx:thmm}
are proved in the supplemental article \cite{Donoho2013c}. The supplemental
article also contains
a derivation of the Stein unbiased risk
estimate for SVST, which is instrumental in the proof of Theorem~\ref{lf:thmm},
and other technical auxiliary lemmas.\vadjust{\goodbreak}

\section{Preliminaries} \label{prelim:sec}

\subsection{Scaling}

Our main object of interest, the worst-case MSE of SVST,
%
\begin{equation}
\label{risk:eq} \mathop{\sup_{X_0\in\Mmn}}_{\operatorname{rank}(X_0)\leq\rho\m}
\frac{1}{\m\n}\E\Vert\hat{X}-X_0 \Vert_F^2
,
\end{equation}
is more conveniently expressed using a specially calibrated risk function.
Since the SVST denoisers are scale-invariant, namely
\[
\E_{X}\bigl\Vert X-\hat{X}(X+\sigma Z)\bigr \Vert^2_F
= \sigma^2 \E_{X}\biggl\Vert\frac{X}{\sigma}-\hat{X} \biggl(
\frac{X}{\sigma
}+Z \biggr) \biggr\Vert^2_F ,
\]
we are free to introduce the scaling
$\sigma=n^{-1/2}$ and define the risk function of a denoiser
$\hat{X}\dvtx \Mmn\to\Mmn$ at
$X_0\in\Mmn$ by
%
\begin{equation}
R(\hat{X},X_0) := \frac{1}{\m}\E\biggl\Vert\hat{X}
\biggl(X_0 + \frac
{1}{\sqrt{\n}}Z \biggr) - X_0
\biggr\Vert_F^2. \label{risk:eqn}
\end{equation}
Then, the worst-case MSE of $\hat{X}$ at $X_0$ is given by
%
\begin{equation}
\label{risk:eq} \mathop{\sup_{X_0\in\Mmn}}_{\operatorname{rank}(X_0)\leq\rho\m}
\frac{1}{\m\n}\E\Vert\hat{X}-X_0 \Vert_F^2
= \mathop{\sup_{{X_0\in\Mmn}}}_{\operatorname{rank}(X_0)\leq\rho\m} R(\hat{X},X_0).
\end{equation}

To vary the SNR in the problem, it will be convenient to vary the norm
of the
signal matrix $X_0$ instead, namely, to consider
$Y=\mu X_0+\frac{1}{\sqrt{\n}}Z$ with $\frac{1}{\m}\Vert X_0 \Vert
_F^2 =1$.

\subsection{Notation}
Vectors are denoted by boldface lowercase letters, such as $\mathbf{v}$,
and their
entries by $v_i$. Matrices are denoted by uppercase letters, such as
$A$, and
their entries by $A_{i,j}$.
Throughout this text, $Y$ will denote the data matrix $Y= \mu
X_0+\frac{1}{\sqrt{n}}Z$. We use $\Mmn$ and $O_\m$ to denote the
set of real-valued
$\m$-by-$\n$ matrices, and group of $\m$-by-$\m$ orthogonal matrices,
respectively. $\Vert\cdot \Vert_F$ denotes the Frobenius matrix norm on
$\Mmn$,
namely the Euclidean norm of a matrix considered as a vector in $\R
^{\m\n}$.
We denote matrix multiplication by either $AB$ or $A\cdot B$.
We use the following convenient notation for matrix diagonals: for a matrix
$X\in\Mmn$, we denote by $X_\on\in\R^\m$ its main diagonal,
%
\begin{equation}
(X_\on)_i = X_{i,i},\qquad 1\leq i \leq\m.
\end{equation}
Similarly, for a vector
$\mathbf{x}\in\R^\m$, and $\n\geq\m$ that we suppress in our
notation, we denote by
$\mathbf{x}_\on\in\Mmn$ the ``diagonal'' matrix
%
\begin{equation}
(\mathbf{x}_\on)_{i,j} = \cases{ x_i, &\quad $1\leq
i=j \leq\m,$ \vspace*{2pt}
\cr
0, &\quad $\mbox{otherwise}.$}
\end{equation}

We use a ``fat'' singular value decomposition (SVD) of $X\in\Mmn$
$X=U_X\cdot\mathbf{x}_\on\cdot V_X'$, with $U_X\in\Mmm$ and
$V_X\in
\Mnn$. Note
that the SVD is not uniquely determined, and in particular $\mathbf
{x}$ can contain
the singular values of $X$ in any order. Unless otherwise noted, we will
assume that the entries of $\mathbf{x}$ are nonnegative and sorted in
nonincreasing order, $x_1\geq\cdots\geq x_\m\geq0$. When $\m< \n
$, the
last $\n-\m$ columns of $V_Y$ are not uniquely determined; we will
see that
our various results do not depend on this choice. Note that
with the ``fat'' SVD, the matrices $Y$ and $U_Y'\cdot Y \cdot V_Y$ have the
same dimensionality, which simplifies the notation we will need.

When appropriate, we let univariate functions act on vectors entry-wise,
namely, for $\mathbf{x}\in\R^\n$ and $f\dvtx \R\to\R$, we write
$f(\mathbf{x})\in\R^\n$ for
the vector with entries $f(\mathbf{x})_i = f(x_i)$.

\subsection{\texorpdfstring{$\hat{X}_\lambda$}{hat{X}lambda} acts by soft thresholding of the data singular
values}

By orthogonal invariance of the Frobenius norm,
\eqref{svt:eq} is equivalent to
%
\begin{equation}
\label{SVST-singvals:eq} \hat{\mathbf{x}}_\lambda=
\mathop{\operatorname{argmin}}_{\mathbf{x}\in\R^\n}
\tfrac{1}{2}\Vert\mathbf{y}-\mathbf{x} \Vert_2^2 +
\lambda\Vert \mathbf{x} \Vert_1 ,
\end{equation}
through the relation $\hat{X}_\lambda(Y) = U_Y \cdot(\hat{\mathbf
{x}}_\lambda)_\on\cdot V_Y'$.
It is well known that the solution to \eqref{SVST-singvals:eq} is
given by
$\hat{\mathbf{x}}_\lambda= \mathbf{y}_\lambda$,
where
$\mathbf{y}_\lambda= (\mathbf{y}-\lambda)_+$ denotes
coordinate-wise soft
thresholding of $\mathbf{y}$ with threshold $\lambda$.
The SVST estimator \eqref{svt:eq} is therefore given by~\cite{Koltchinskii2011a}
%
\begin{equation}
\label{svt2:eq} \hat{X}_\lambda\dvtx Y\mapsto U_Y \cdot(
\mathbf{y}_\lambda)_\on\cdot V_Y'.
\end{equation}

Note that \eqref{svt2:eq} is well defined, that is, $\hat{X}_\lambda
(Y)$ does
not depend on the particular SVD $Y=U_Y\cdot(\mathbf{y})_\on\cdot
V_Y'$ chosen.

In case $\mathrm{Sym}$, observe that the solution to \eqref{svt:eq} is
constrained to lie
in the linear subspace of symmetric matrices. The solution is the same whether
the noise matrix $Z\in\Mnn$ has i.i.d. standard normal entries, or
whether $Z$ is
a symmetric Wigner matrix $ \frac{1}{2}(Z_1+Z_1')$ where $Z_1\in\Mnn
$ has
i.i.d. standard normal entries. Below, we assume that the data in case
$\mathrm{Sym}$ is
of the form $Y=\mu X_0+Z/
\sqrt{\n}$ where $X_0\in S_+^\n$ and $Z$ has this Wigner form, namely,
the singular values $\mathbf{y}$ are the absolute values of
eigenvalues of the
symmetric matrix $Y$.

\section{The least-favorable matrix for SVST is at \texorpdfstring{$\Vert X\Vert=\infty$}
{$||X||=infty$}} \label{lf:sec}

We now prove Theorem~\ref{lf:thmm}, which characterizes the worst-case
MSE of the
SVST denoiser $\hat{X}_\lambda$ for a given $\lambda$. The theorem
follows from
a combination of two classical gems of the statistical literature. The
first is
Stein's unbiased risk estimate (SURE) from 1981, which we
specialize to the SVST estimator; see also \cite{Cand2012}. The second is
Anderson's celebrated monotonicity property for the integral of a symmetric
unimodal probability distribution over a symmetric convex set
\cite{Anderson1955}, from 1955, and more specifically its implications for
monotonicity of the power function of certain tests in multivariate hypothesis
testing \cite{Das1964}.
To simplify the proof, we introduce the following definitions, which
will be
used in this section only.

\begin{defn}[(A weak notion of matrix majorization based on singular values)] \label{majorization:def}
Let $A,B\in\Mmn$ have singular value vectors $\mathbf{a},\mathbf
{b}\in\R^\m$,
respectively, which as usual we assume are sorted in nonincreasing order:
$0\leq a_\m\leq\cdots\leq a_1$ and
$0\leq b_\m\leq\cdots\leq b_1$.
If $a_i \leq b_i$ for $i=1,\ldots,\m$, we write $A\preceq B$.
\end{defn}

We note that by rescaling an arbitrary rank-$\rk$
matrix, it is always possible to majorize any fixed matrix of rank at
most $\rk$ (in the sense of
Definition~\ref{majorization:def}).

\begin{lemma} \label{mu:lem}
Let $C\in\Mmn$ be a matrix of rank $\rk$,
and let $X\in\Mmn$ be a matrix of
rank at most $\rk$. Then there exists
$\mu>0$ for which $X\preceq\mu C$.
\end{lemma}

\begin{pf}
Let $\mathbf{c},\mathbf{x}$ be the vectors of singular values of
$C,X$, respectively,
each sorted in
nonincreasing order.
Then $c_\rk>0$. Take $\mu= x_1/c_\rk$.
For $1\leq i \leq\rk$ we have $x_i\leq x_1 = \mu c_\rk\leq\mu c_i$,
and for
$\rk+1\leq i\leq\m$ we have $\mu c_i=x_i=0$.
\end{pf}

The above weak notion of majorization gives rise to a weak notion of
monotonicity:

\begin{defn}[(Orthogonally invariant function of a matrix argument)]
We say that $f\dvtx \Mmn\to\R$ is an orthogonally invariant function if
$f(U\cdot A\cdot V')=f(A)$ for all $A\in\Mmn$ and all orthogonal
$U\in O_\m$ and $V\in O_\n$.
\end{defn}

\begin{defn}[(SV-monotone increasing function of a matrix argument)]
Let $f\dvtx \Mmn\to\R$ be orthogonally invariant.
If, whenever $A\preceq B$ and $\sigma>0$, $f$ satisfies
%
\begin{equation}
\label{matrix-monotone:eq} \E f(A+Z) \leq\E f(B+Z) ,
\end{equation}
for $Z\in\Mmn$ and $Z_{i,j}\stackrel{\mathrm{i.i.d.}}{\sim}\Nc
(0,\sigma^2)$, we say that $f$ is
singular-value-monotone increasing, or \textit{SV-monotone increasing}.
\end{defn}

We now provide a sufficient condition for SV-monotonicity, which
follows from
Anderson's seminal monotonicity result
\cite{Anderson1955}. The following lemma is proved in the supplemental article
\cite{Donoho2013c}.

\begin{lemma} \label{sufficient-sv-mon-II:lem}
Assume that $f\dvtx \Mmn\to\R$ can be decomposed as $f=\sum_{k=1}^s
f_k$, where
for each $1\leq k\leq s$, $f_k\dvtx \Mmn\to\R$ is
a bounded, orthogonally invariant function. Further assume that for each
$1\leq k \leq s$, $f_k$ is
quasi-convex, in the sense that for all $c\in\R$, the set
$f_k^{-1}((-\infty,c])$ is convex in $\Mmn$.
Then $f$ is SV-monotone increasing.
\end{lemma}

The second key ingredient in the proof of Theorem~\ref{lf:thmm} is the Stein
unbiased risk estimate for SVST.
Let $\hat{X}$ be a weakly differentiable estimator of $X_0$ from data
$Y=X_0 +
\sigma Z$, where $Z$ has i.i.d. standard normal entries. The Stein
unbiased risk
estimate \cite{Stein1981} is a function of the data, $Y\mapsto
\operatorname{SURE}(Y)$, for
which $\E \operatorname{SURE}(Y) = \E\Vert\hat{X}-X_0 \Vert_2^2$. In our case, $X_0,
Z$ and $Y$
are matrices in $\Mmn$, and Stein's theorem (\cite{Stein1981}, Theorem~1) implies
that for
\[
\operatorname{SURE}(Y) = \m\n\sigma^2 + \bigl\Vert\hat{X}(Y) - Y\bigr \Vert^2_F
+ 2\sigma^2\sum_{i,j}\frac{\partial( \hat
{X}(Y) -
Y)_{i,j}}{\partial Y_{i,j}} ,
\nonumber
\]
we have
\[
\Vert\hat{X}-X_0 \Vert_F^2 =
\E_{X_0} \operatorname{SURE}(Y).
\]

In the supplemental article \cite{Donoho2013c}, we derive SURE for a
large class
of invariant matrix denoisers. As a result, we prove:

\begin{lemma}[(The Stein unbiased risk estimate for SVST)] \label{svt-SURE:lem}
For each $\lambda>0$, there exists an event $\mathcal{S}\subset\Mmn$
and a function, 
$\operatorname{SURE}_\lambda\dvtx\mathcal{S}\to\R$ which maps a matrix $Y$ with
singular values
$\mathbf{y}$ to
\begin{eqnarray*}
\operatorname{SURE}_\lambda(Y) &=& \m
\nonumber
+\sum_{i=1}^\m
\biggl[ \bigl(\min \{ y_i,\lambda \}\bigr)^2 -
\textbf{1}_{\{y_i<\lambda\}} - \frac{(\n-\m)\cdot\min \{ y_i,\lambda \}}{y_i} \biggr] \label{SURE-svt-rest:eq}
\\
&&{}-\frac{2}{\n} \sum_{1\leq i\neq j\leq\m} \frac{\min \{ y_j,\lambda
 \}y_j-\min \{ y_i,\lambda
 \}y_j}{y_j^2-y_i^2}
, \label{SURE-svt-div:eq}
\end{eqnarray*}
enjoying the following properties:
\begin{longlist}[(1)]
\item[(1)] $\mathbb{P}(\mathcal{S})=1$, where $\mathbb{P}$ is the
distribution of the
matrix $Z$ with $Z_{i,j}\stackrel{{i.i.d.}}{\sim}\Nc(0,1)$.

\item[(2)] $\operatorname{SURE}_\lambda$ is a finite sum of bounded, orthogonally
invariant, quasi-convex functions.
\item[(3)] Denoting as usual $Y=X_0+Z/\sqrt{\n}\in\Mmn$, where
$X_0,Z\in\Mmn$ and $Z_{i,j}\stackrel{{i.i.d.}}{\sim}\Nc
(0,1)$, we have
\[
R(\hat{X}_\lambda,X_0) = \frac{1}{\m}
\E_{X_0} \operatorname{SURE}_\lambda(Y).
\]
\end{longlist}
\end{lemma}

Putting together Lemmas \ref{sufficient-sv-mon-II:lem} and~\ref{svt-SURE:lem},
we come to a crucial property of SVST.

\begin{lemma}[(The risk of SVST is monotone nondecreasing in the signal singular
values)] \label{svt-mon:lem}
For each $\lambda>0$, the map $X\mapsto R(\hat{X}_\lambda,X)$ is a bounded,
SV-monotone increasing function. In particular, let
$A,B\in\Mmn$ with $A \preceq B$. Then
%
\begin{equation}
R(\hat{X}_\lambda,A) \leq R(\hat{X}_\lambda,B).
\end{equation}
\end{lemma}

\begin{pf}
By Lemma~\ref{svt-SURE:lem}, the function $\operatorname{SURE}_\lambda\dvtx \Mmn\to\R$
satisfies the conditions of Lemma~\ref{sufficient-sv-mon-II:lem} and is
therefore SV-monotone increasing. It follows that
\begin{eqnarray*}
R(\hat{X}_\lambda,A) &=& \frac{1}{\m}\E_A
\operatorname{SURE}_\lambda(A+Z/\sqrt{\n}) \\
&\leq& \frac{1}{\m}\E_B
\operatorname{SURE}_\lambda(B+Z/\sqrt{\n})= R(\hat{X}_\lambda,B).
\end{eqnarray*}
To see that the risk is bounded, note that for any $X\in\Mmn$, we
have by Lemma~\ref{svt-SURE:lem}
\[
\infty<\inf_{Y\in\Mmn} \E \operatorname{SURE}_\lambda(Y) \leq R(
\hat{X}_\lambda,X) \leq\sup_{Y\in\Mmn} \E
\operatorname{SURE}_\lambda(Y) < \infty.
\]
\upqed\end{pf}

\begin{pf*}{Proof of Theorem \protect\ref{lf:thmm}}
By Lemma~\ref{svt-mon:lem}, the map
$\mu\to R(\hat{X}_\lambda,\mu C)$ is bounded and monotone
nondecreasing in
$\mu$. Hence $\lim_{\mu\to\infty} R(\hat{X}_\lambda,\mu C)$
exists and is
finite, and
%
\begin{equation}
\label{lim-larger:eq} R(\hat{X}_\lambda,\mu_0 C) \leq\lim
_{\mu\to\infty} R(\hat{X}_\lambda,\mu C)
\end{equation}
for all
$\mu_0>0$. Since $\operatorname{rank}(C) = r$, obviously
\[
\sup_{
\operatorname{rank}(X_0)\leq\rk} R(\hat{X}_\lambda,X_0) \geq
\lim_{\mu\to\infty} R(\hat{X}_\lambda, \mu C) ,
\]
and we only need to show the
reverse inequality. Let $X_0\in\Mmn$ be an arbitrary matrix of rank
at most
$\rk$.
By Lemma~\ref{mu:lem} there exists $\mu_0$ such that $X_0\preceq\mu
_0 C$. It
now follows from Lemma~\ref{svt-mon:lem} and \eqref{lim-larger:eq} that
\[
R(\hat{X}_\lambda,X_0) \leq R(\hat{X}_\lambda,\mu
_0 C) \leq \lim_{\mu\to\infty} R(\hat{X}_\lambda,
\mu C).
\]
\upqed\end{pf*}

\section{Worst-case MSE} \label{finite-n-proxy:sec}

Let $\lambda$ and $\rk\leq\m\leq\n$, and consider them fixed for
the remainder
of this section. Our second main result,
Theorem~\ref{finite-n-mmx:thmm}, follows immediately
from Theorem~\ref{lf:thmm}, combined with the
following lemma, which is proved in the supplemental article \cite{Donoho2013c}.

\begin{lemma} \label{proxy:lem}
Let $X_0\in\Mmn$ be of rank $\rk$.
Then
\[
\lim_{\mu\to\infty} R(\hat{X}_\lambda, \mu X_0) =
\MSE_{\n} \biggl(\frac{\lambda}{\sqrt{1-\rk/\n}} ; \rk,\m ,\alpha \biggr) ,
\]
as defined in \eqref{finite-n-proxy:eq},
with $\alpha=1$ for case $\mathrm{Mat}$ and $\alpha=1/2$ for case $\mathrm{Sym}$.
\end{lemma}

In the supplemental article \cite{Donoho2013c} we prove the following lemma:

\begin{lemma} \label{finite-n-MSE-convex:lem}
The function $\Lambda\mapsto\MSE_\n(\Lambda;\rk,\m,\alpha)$,
defined in
\eqref{finite-n-proxy:eq}
on
$\Lambda\in[0,\infty)$, is convex and obtains a unique minimum.
\end{lemma}

Our second main result is an immediate consequence:

\begin{pf*}{Proof of Theorem \protect\ref{finite-n-mmx:thmm}}
Let $C\in\Mmn$ be an arbitrary fixed matrix of rank~$\rk$.
For case $\mathrm{Mat}$, by Theorem~\ref{lf:thmm} and Lemma~\ref{proxy:lem},
\begin{eqnarray*}
\mmx_\n(\rk,\m|\mathrm{Mat})&=& \inf_\lambda\mathop{\sup
_{X_0\in\Mmn}}_{\operatorname{rank}(X_0)\leq\rk} R(\hat{X}_\lambda,X_0)
= \inf_{\lambda>0} \lim_{\mu\to\infty} R(
\hat{X}_\lambda, \mu C)
\\
&=& \inf_{\lambda>0} \MSE_{\n} \biggl(\frac{\lambda}{\sqrt{1-\rk/\n}} ;
\rk,\m,1 \biggr) \\
&=& \min_{\Lambda>0} \MSE_\n(\Lambda;\rk,
\m,1) ,
\end{eqnarray*}
where we have used Lemma~\ref{finite-n-MSE-convex:lem}, which also asserts
that the minimum is unique.

Now let $C\in S_+^\n$ be an arbitrary, fixed symmetric positive semidefinite
matrix of rank $\rk$.
For case $\mathrm{Sym}$, by the same lemmas,
\begin{eqnarray*}
\mmx_n(\rk|\mathrm{Sym})&=& \inf_\lambda\mathop{\sup
_{X_0\in\Mmn}}_{\operatorname{rank}(X_0)\leq\rk} R(\hat{X}_\lambda,X_0)
= \inf_\lambda \lim_{\mu\to\infty} R(
\hat{X}_\lambda, \mu C)
\\
&=& \inf_\lambda \MSE_{\n} \biggl(\frac{\lambda}{\sqrt{1-\rk/\n}}
; \rk,1/2 \biggr) = \min_\Lambda\MSE_\n(\Lambda;
\rk,1/2).
\end{eqnarray*}
\upqed\end{pf*}

\section{Worst-case AMSE} \label{asymp-proxy:sec}
Toward the proof of our third main result, Theorem~\ref
{asymp-mmx:thmm}, let
$\lambda$ be fixed. We first show that in the proportional growth framework,
where the rank $\rk(\n)$, number of rows $\m(\n)$ and
number of columns $\n$ all tend to $\infty$ proportionally to each other,
the key quantity in our formulas can be evaluated by complementary incomplete
moments of a Mar\u{c}enko--Pastur distribution, instead of a sum of
complementary incomplete moments of Wishart eigenvalues.

\begin{defn} \label{zeta:def}
For a pair of matrices $X_0,Z\in\Mmn$, we denote by $\zeta
(X_0, Z|\mathrm{Mat}) =
(\zeta_1 ,\ldots, \zeta_{\m-\rk})$ the singular values, in nonincreasing
order, of
%
\begin{equation}
\Pi_\m\cdot Z \cdot\Pi_\n'\in
M_{(\m-\rk)\times(\n-\rk)} ,
\end{equation}
where $\Pi_\m\dvtx \R^\m\to\R^{\m-\rk}$ is the projection of $\R^\m
$ on
$\operatorname{null}(X_0')=\operatorname{Im}(X_0)^\perp$ and $\Pi_\n\dvtx \R^\n\to\R^{\n-\rk}$ is
the projection on
$\operatorname{null}(X_0)$. Similarly, for a pair of matrices $X_0,Z\in\Mnn$, denote by
$\zeta(X_0,Z|\mathrm{Sym}) = (\zeta_1 ,\ldots, \zeta_{\m-\rk})$ the
eigenvalues, in
nonincreasing order, of
%
\begin{equation}
\Pi_\m\cdot\tfrac{1}{2}\bigl(Z+Z'\bigr) \cdot
\Pi_\n'\in M_{(\n-\rk)\times(\n-\rk)}.
\end{equation}
\end{defn}

\begin{lemma} \label{MP-asymp:lem}
Consider sequences $\n\mapsto\rk(\n)$ and $\n\mapsto\m(\n)$ and numbers
$0<\beta\leq1$ and $0\leq\rho\leq1$ such that
$\lim_{\n\to\infty}\rk(\n)/\m(\n) = \rho$
and $\lim_{\n\to\infty}\m(\n)/\n= \beta$.
Let $(\zeta_1(\n),\ldots,\zeta_{\m-\rk}(\n))=\zeta(X_0,Z|\bX
)$, as in Definition~\ref{zeta:def}, where $Z\in\Mmn$ has i.i.d. $\Nc(0,1)$
entries.
Define $\gamma=(\beta-\rho\beta)/(1-\rho\beta)$ and
$\gamma_{\pm}= ( 1\pm\sqrt{\gamma}  )^2$, and let
$0\leq\Lambda\leq\sqrt{\gamma_+}$.
Then
\[
\lim_{\n\to\infty} \frac{1}{\m}\sum
_{i=1}^{\m-\rk} \E \biggl(\frac{\z_i}{\sqrt{\n-\rk}} - \Lambda
\biggr)^2_+ = (1-\rho)\int_{\Lambda^2}^{\gamma_+}(\sqrt{t}-
\Lambda)^2 \frac{\sqrt{(\gamma_+ -t)(t-\gamma_-)}}{2\pi t\gamma} \,dt.
\]
\end{lemma}
\begin{pf}
Write $\xi_i = \zeta^2_i /(\n-\rk)$, and recall that
by the Mar\u{c}enko--Pastur law~\cite{Marcenko1967},
\[
\lim_{\n\to\infty} \frac{1}{\m-\rk}\sum
_{i=1}^{\m-\rk} \delta _{\xi_i} \stackrel{w} {=}
P_\gamma,
\]
in the sense of weak convergence of probability measures,
where $P_\gamma$ is the Mar\u{c}enko--Pastur probability
distribution
with density $p_\gamma=dP_\gamma/dt$ given by~\eqref{mp:eq}.
Now,
\begin{eqnarray*}
\lim_{\n\to\infty} \frac{1}{\m} \sum
_{i=1}^{\m-\rk} (\sqrt{\xi_i} - \Lambda
)^2_+ &=& \lim_{\n\to\infty} \frac{1}{\m} \sum
_{i=1}^{\m-\rk}\int_0^\infty
(\sqrt{t} - \Lambda )^2_+ \delta_{\xi_i}(t) \,dt
\\
&=& \lim_{\n\to\infty} \biggl(1-\frac{\rk}{\m} \biggr) \int
_0^\infty (\sqrt{t} - \Lambda )^2_+
\frac{1}{\m-\rk} \sum_{i=1}^{\m-\rk}
\delta_{\xi_i}(t) \,dt
\\
&=& (1-\rho ) \int_0^{\gamma_+} (\sqrt{t} - \Lambda
)^2_+ p_\gamma(t)\,dt
\end{eqnarray*}
as required.
\end{pf}

\begin{lemma}
\label{asymp-proxy:lem}
Let $\m(\n)$ and $\rk(\n)$ such that
$\lim_{\n\to\infty} \m/\n=\beta$ and\break $\lim_{\n\to\infty} \rk
/\m=\rho$,
and set $\tilde{\rho} = \beta\rho$.
Then
\[
\lim_{\n\to\infty}\mathop{\sup_{
X_0\in\Mmn}}_{{\operatorname{rank}(X_0)\leq\rk}}
R(\hat{X}_\lambda,X_0) = \MSE \biggl(\frac{\lambda}{\sqrt{1-\tilde{\rho}}};
\rho,\tilde {\rho},\alpha \biggr) ,
\]
where the right-hand side is defined in \eqref{proxy:eq}, with $\alpha
=1$ for
case $\mathrm{Mat}$ and $\alpha=1/2$ for case $\mathrm{Sym}$.
\end{lemma}

\begin{pf}
For case $\mathrm{Mat}$, let $C(\n)\in\Mmn$ be an arbitrary fixed matrix of
rank $\rk$.
For case $\mathrm{Sym}$,
$C(\n)\in S_+^\n$ an arbitrary, fixed symmetric positive semidefinite
matrix of rank $\rk$.
By Theorem~\ref{lf:thmm} and Lemma~\ref{proxy:lem},
\begin{eqnarray*}
&& \lim_{\n\to\infty}\mathop{\sup_{X_0\in\Mmn}}_{\operatorname{rank}(X_0)\leq\rk}
R(\hat{X}_\lambda,X_0)\\
&&\qquad= \lim_{\n\to\infty} \lim
_{\mu\to\infty} R\bigl(\hat{X}_\lambda, \mu C(\n)\bigr)
\\
&&\qquad= \lim_{\n\to\infty} \Biggl[ \frac{\rk}{\m}+\frac{\rk}{\n}-
\frac{r^2}{\m\n} + \frac{\rk}{\m}\lambda^2\\
&&\hspace*{60pt}{} + \alpha
\frac{\n-\rk}{\m\n}\sum_{i=1}^{\m-\rk} \E \biggl(
\frac{\zeta_i}{\sqrt{\n-\rk}} - \frac{\lambda}{\sqrt{1-\rk/\n}}
 \biggr)^2_+ \Biggr]
\\
&&\qquad= \rho+ \tilde{\rho} - \rho\tilde{\rho} + (1-\tilde{\rho}) \rho
\Lambda^2 \\
&&\qquad\quad{}+ \alpha(1-\rho) (1-\tilde{\rho}) \int_{\Lambda^2}^{\gamma_+}
(\sqrt{t}-\Lambda)^2 MP_\gamma(t)\,dt
\\
&&\qquad=\MSE \biggl(\frac{\lambda}{\sqrt{1-\tilde{\rho}}};\rho,
\tilde {\rho},\alpha \biggr),
\end{eqnarray*}
where we have used Lemma~\ref{MP-asymp:lem} and set
$\Lambda=\lambda/\sqrt{1-\tilde{\rho}}$.
\end{pf}

In the supplemental article we prove a variation of Lemma~\ref{finite-n-MSE-convex:lem} for the asymptotic setting:

\begin{lemma} \label{asymp-MSE-convex:lem}
The function $\Lambda\mapsto\MSE(\Lambda;\rho,\tilde{\rho
},\alpha)$, defined in
\eqref{proxy:eq}
on
$\Lambda\in[0,\gamma_+]$, where
$\gamma_+ =  (1+\sqrt{ (\tilde{\rho} - \rho\tilde{\rho}) /
(\rho-\rho\tilde{\rho})} )^2$,
is convex and obtains a unique minimum.
\end{lemma}

This allows us to the prove our third main result.
\begin{pf*}{Proof of Theorem \protect\ref{asymp-mmx:thmm}}
By Lemma~\ref{asymp-proxy:lem},
\begin{eqnarray*}
\mmx(\rho,\beta|\X) &=& \lim_{\n\to\infty}\inf_\lambda
\mathop{\sup_{X_0\in\Mmn}}_{{\operatorname{rank}(X_0)\leq\rk}} R(\hat{X}_\lambda,X_0)
\\
&=& \inf_\lambda\lim_{\n\to\infty} \mathop{\sup
_{
X_0\in\Mmn}}_{{\operatorname{rank}(X_0)\leq\rk}} R(\hat{X}_\lambda,X_0)
\\
&=& \inf_\lambda \MSE \biggl(\frac{\lambda}{\sqrt{1-\tilde{\rho}}};\rho,\tilde {
\rho},\alpha \biggr) = \min_\Lambda\MSE (\Lambda;\rho,\tilde{
\rho},\alpha ) ,
\end{eqnarray*}
with $\alpha=1$ for case $\mathrm{Mat}$ and $\alpha=1/2$ for case $\mathrm{Sym}$,
where we have used Lemma~\ref{asymp-MSE-convex:lem}, which also
asserts that the
minimum is unique.
\end{pf*}

\section{Minimax AMSE} \label{asymp-minimizer:sec}

Having established that the asymptotic worst-case MSE~\eqref{proxy:eq}
satisfies
\eqref{asymp-mmx-mat:eq} and \eqref{asymp-mmx-sym:eq}, we turn to its
minimizer $\Lambda_*$. The notation follows~\eqref{proxy-minimizer:eq}.

\begin{pf*}{Proof of Theorem \protect\ref{proxy-minimizer:thmm}}
By equation (4.2) in the supplemental article
\cite{Donoho2013c}, the condition
\[
\frac{d\mathbf{M}(\Lambda;\rho,\tilde{\rho},\alpha)}{d\Lambda}=0
\]
is thus equivalent, for any $\rho\in[0,1]$, to
%
\begin{equation}
\label{Lambda:eq} f(\Lambda,\rho) := \rho\Lambda- \alpha(1-\rho)
\int_{\Lambda^2}^{\gamma_+} (\sqrt{t}-\Lambda) p_\gamma(t) \,dt =0
,
\end{equation}
establishing \eqref{argmin-Lambda:eq} in particular for
$0<\rho<1$.
By Lemma~\ref{asymp-MSE-convex:lem}, the minimum exists and is unique; namely
this equation has a unique
root in $\Lambda$.
One directly verifies that
$f(1+\sqrt{\beta},0)=f(0,1)=0$. The limits \eqref
{lim_Lambda_down:eq} and
\eqref{lim_Lambda_up:eq} follow from the fact that
$\rho\mapsto\Lambda_*(\rho,\cdot)$ is decreasing. To establish this,
it is enough to observe that
$\partial f/\partial\rho> 0$ for all $(\Lambda,\rho)$, which can be verified
directly.
\end{pf*}

Theorem~\ref{proxy-sq:thmm}, which provides more a explicit formula
for the
minimax AMSE in square matrix case ($\beta=1$), is proved in the supplemental
article \cite{Donoho2013c}.

\section{Global minimax MSE and AMSE} \label{global:sec}

In this section we prove Theorem~\ref{global-mmx:thmm}, which provides
a lower
bound on the minimax risk of the family of all measurable matrix
denoisers (as opposed
to the family of SVST denoisers considered so far) over $\m$-by-$\n$
matrices of rank at most $\rk$.
Consider the class of singular-value matrix denoisers, namely
all mappings $Y\mapsto\hat{X}(Y)$
that act on the data $Y$ only through their singular values. More specifically,
consider all denoisers
$\hat{X}\dvtx \Mmn\to\Mmn$ of the form
%
\begin{equation}
\label{sv-denoiser:eq} \hat{X}(Y) = U_Y\cdot\hat{\mathbf{x}}(
\mathbf{y})_\on\cdot V_Y' ,
\end{equation}
where $Y=U_Y \cdot\mathbf{y}_\on\cdot V_Y'$ and
$\hat{\mathbf{x}}\dvtx [0,\infty)^\m\to[0,\infty)^\m$.
(Note that this class contains SVST denoisers but does not exhaust all
measurable denoisers.)
The mapping in~\eqref{sv-denoiser:eq} is not well defined in general, since
the SVD of $Y$, and in particular the order of the singular values in the
vector $\mathbf{y}$, is not uniquely determined.
However, \eqref{sv-denoiser:eq}
is well defined when each function $\hat{x}_i\dvtx [0,\infty)\to[0,\infty)$
is invariant under permutations of its coordinates. Since the equality
$Y=U_Y \cdot\mathbf{y}_\on\cdot V_Y'$ may hold for vectors $\mathbf
{y}$ with negative
entries, we are led to the following definition.
%
\begin{defn} \label{singval-denoiser:def}
By \textit{singular-value denoiser} we mean any
measurable mapping
$\hat{X}\dvtx \Mmn\to\Mmn$ which takes the form \eqref{sv-denoiser:eq},
where each entry
of $\hat{\mathbf{x}}$ is a function
$\hat{x}_i\dvtx \R^\m\to\R$ that is invariant under permutation and
sign changes of its
coordinates. We let $\cD$ denote the class of such mappings.
\end{defn}

For a detailed introduction to real-valued or matrix-valued functions which
depend on a matrix argument only through its singular values, see \cite{Lewis1995,Lewis2005}.
The following lemma is proved in the supplemental article \cite{Donoho2013c}.

\begin{lemma}[(Singular-value denoisers can only improve in worst-case)] \label{huntstein:lem}
Let $\hat{X}_1\dvtx \Mmn\to\Mmn$ be an arbitrary measurable matrix
denoiser. There
exists a singular-value denoiser $\hat{X}$ such that
\[
\mathop{\sup_{{X_0\in\Mmn}}}_{{\operatorname{rank}(X_0)\leq\rk}} R(\hat{X},X_0)
\leq \mathop{\sup_{X_0\in\Mmn}}_{\operatorname{rank}(X_0)\leq\rk} R(\hat{X}_1,X_0).
\]
\end{lemma}

\begin{pf*}{Proof of Theorem \protect\ref{global-mmx:thmm}}
We consider the case $\bX=\mathrm{Mat}_{\m,\n}$.
By Lemma~\ref{huntstein:lem}, it is enough to show that
\[
\frac{\rk}{\m}+\frac{\rk}{\n}-\frac{\rk^2+\rk}{\m\n} \leq \mathop{\sup
_{X_0\in\Xmn}}_{{\operatorname{rank}(X_0)\leq\rk}} R(\hat{X},X_0) ,
\]
where $\hat{X} \in\cD$ is an arbitrary singular-value denoiser.
Indeed, let $X_0\in\Mmn$ be a fixed arbitrary matrix of rank $\rk$.
The calculation leading to equation (3.9) in the
supplemental article
\cite{Donoho2013c} is valid for any
rule in $\cD$, and implies that
$R(\hat{X}(Y),X_0) \geq1 - \frac{1}{\m}\E\Vert\mathbf{z} \Vert_2^2$,
where $Y=U_Y\cdot\mathbf{y}_\on\cdot V_Y'$ and
%
\begin{equation}
\mathbf{z} = \frac{1}{\sqrt{\n}}\bigl(U_Y' \cdot Z
\cdot V\bigr)_\on.
\end{equation}
Write $Y_\mu= \mu X_0 + Z/\sqrt{\n}=U_\mu\cdot(\mathbf{y}_\mu
)_\on
\cdot V_\mu'$,
and let $\mathbf{z}_\mu= \frac{1}{\sqrt{\n}}(U_\mu' \cdot Z \cdot
V_\mu)_\on$. We therefore have
\[
\mathop{\sup_{X_0\in\Xmn}}_{\operatorname{rank}(X_0)\leq\rk} R(\hat{X},X_0)
\geq \lim_{\mu\to\infty} R(\hat{X},\mu X_0) \geq 1 -
\frac{1}{\m} \lim_{\mu\to\infty} \E\Vert z_\mu
\Vert_2^2.
\]

Combining equations (3.17) and (3.15) in the
supplemental article \cite{Donoho2013c}, we have
\[
\frac{1}{\m} \sum_{i=\rk+1}^\m\lim
_{\mu\to\infty} \E(z_{\mu
,i})^2 = 1 -
\frac{\rk}{\m}-\frac{\rk}{\n} + \frac{\rk^2}{\m\n}.
\]
A similar argument yields
$ \frac{1}{\m}\sum_{i=1}^\rk\lim_{\mu\to\infty} \E(z_{\mu
,i})^2 =
\frac{\rk}{\m\n}$,
and the first part of the theorem follows. The second part of the theorem
follows since, taking the limit $\n\to\infty$ as prescribed, we have
$\rk/\m\to\rho$, $\rk/\n\to\tilde{\rho}$ and $\rk/\m\n\to0$.
For the third part of the theorem, we have by Theorem~\ref{small-rho-mmx:thmm},
\begin{eqnarray*}
\lim_{\rho\to0}\frac{\mmx(\rho,\beta|\bX)}{\mmx^-(\rho,\beta
)} &=& \lim_{\rho\to0}
\frac{\mmx(\rho,\beta|\bX)}{\rho+\beta\rho
+\beta\rho^2} = \frac{2(1+\sqrt{\beta}+\beta)}{1+\beta}\\
& = &2 \biggl( 1+\frac{\sqrt{\beta}}{1+\beta} \biggr).
\end{eqnarray*}
\upqed\end{pf*}

\section{Discussion}

In the \hyperref[intro:sec]{Introduction}, we pointed out several ways that
these matrix denoising results for SVST estimation
of low-rank matrices parallel results for soft thresholding
of sparse vectors. Our derivation of the minimax MSE formulas
exposed two more parallels:
\begin{itemize}
\item\textit{Common structure of minimax MSE formulas.}
The minimax MSE formula vector denoising problem involves certain
incomplete moments of the standard Gaussian distribution \cite{Donoho2011}.
The matrix denoising problem involves completely analogous
incomplete moments, only replacing the Gaussian by the
Mar\v{c}enko--Pastur distribution or (in the square case $\beta=1$)
the quarter-circle law.
\item\textit{Monotonicity of SURE}. In both settings, the
least-favorable estimand places the signal
``at $\infty$,'' which yields a convenient formula for
Minimax MSE \cite{Donoho2011}.
In each setting, validation of the least-favorable estimation
flows from monotonicity, in an appropriate sense, of Stein's unbiased
risk estimate within that specific setting.
\end{itemize}

\section*{Acknowledgments}

We thank Iain Johnstone, Andrea Montanari and Art Owen for advice at
several crucial
points, and the anonymous referees for many helpful suggestions.


\begin{supplement}[id=suppA]
\stitle{Proofs and additional discussion}
\slink[doi]{10.1214/14-AOS1257SUPP} 
\sdatatype{.pdf}
\sfilename{aos1257\_supp.pdf}
\sdescription{In this supplementary material we prove Theorems \ref
{proxy-sq:thmm},
\ref{asymp-tuning-mmx:thmm},
\ref{parametric:thmm},
\ref{small-rho-mmx:thmm} and other lemmas. We also derive the Stein
unbiased risk
Estimate (SURE) for SVST, which is instrumental in the proof of
Theorem~\ref{lf:thmm}. Finally, we discuss similarities
between singular value
thresholding and sparse vector thresholding.}
\end{supplement}

%

\printaddresses
\end{document}